\documentclass{article}
\usepackage{amsmath}
\usepackage{amsfonts}
\usepackage{amsthm}
\usepackage{amssymb}
\usepackage[english]{babel}
\usepackage[ansinew]{inputenc}
\usepackage[all]{xy}
\usepackage{color}

\theoremstyle{definition}     

\newtheorem{defi}{Definition}[section]
\newtheorem{remark}[defi]{Remark}
\newtheorem{notation}[defi]{Notation}
\newtheorem{example}[defi]{Example}
\newtheorem{definition}[defi]{Definition}

\theoremstyle{plain}

\newtheorem{theorem}[defi]{Theorem}
\newtheorem{corollary}[defi]{Corollary}
\newtheorem{lemma}[defi]{Lemma}
\newtheorem{proposition}[defi]{Proposition}

\input xy
\xyoption{all}

\title{Pushout of quasi-finite and flat group schemes over a Dedekind ring}

\author{Marco Antei\thanks{The author was supported by a postdoctoral fellowship funded by the Skirball Foundation via the Center for Advanced Studies in Mathematics at Ben-Gurion University of the Negev. The author was also supported by the Israel Science Foundation (grant No.\ 23/09).}}

\date{}
\begin{document}
\maketitle

\textbf{Abstract.}  Let $G$, $G_1$ and $G_2$ be quasi-finite and flat group schemes over a complete discrete valuation ring $R$, $\varphi_1:G\to G_1$ any morphism of $R$-group schemes and $\varphi_2:G\to G_2$ a model map. We construct the pushout $P$ of  $G_1$ and $G_2$ over $G$ in the category of $R$-affine group schemes. In particular when $\varphi_1$ is a model map too we show that $P$ is still a model of the generic fibre of $G$. We also provide a short proof for the existence of cokernels and quotients of finite and flat group schemes over any Dedekind ring.



\tableofcontents
\bigskip

\section{Introduction}
\subsection{Aim and scope}
We are interested in the
construction of the pushout (whose definition will be recalled in \S \ref{sez:poha}) in the category of affine group schemes over a given ring as described hereafter. It is known that in the category of abstract groups
the pushout of two groups over a third one always exists but it is not finite even
when the three groups are all finite (unless one takes very particular cases). However for group schemes over a Dedekind ring $R$ something new  happens when we consider some special important cases: so let $G$, $G_1$ and $G_2$ be $R$-affine group schemes and consider the diagram 

\begin{equation}\label{eqIntro}
\xymatrix{ & G\ar[ld]_{\varphi_1}\ar[rd]^{\varphi_2} & \\G_1 & & G_2 }
\end{equation} where  $\varphi_i:G\to G_i$ ($i=1,2$) are  $R$-group scheme morphisms. We first prove the following

%As the author was not able to find a complete treatment on the topic (but see for instance \cite{FD}, I, \S 3 for a particular case which inspired our description), the tensor product of non commutative algebras will be recalled in the  Appendix \ref{sez:AppA} 
%mainly to state the universal property of the tensor product of noncommutative algebras (cf. Proposition \ref{propUnivTensor}). It differs from the usual one  for commutative algebras and in particular in the noncommutative case the tensor product is not a pushout (cf. for instance Example \ref{exaCopro2} but one can find easier examples). The diagrammatic definitions of algebra and Hopf algebra will also be recalled in the Appendix as strongly used in this paper.

\begin{theorem} (Cf. Theorem \ref{teoFinito})
Assume $R$ is a complete discrete valuation ring and $G$, $G_1$, $G_2$ are finite and flat over $R$. Then if $\varphi_1$ is a model map (i.e. generically an isomorphism) the pushout of (\ref{eqIntro}) in the category of affine $R$-group schemes exists. Moreover it is finite and flat and its generic fibre is isomorphic to $G_{2,K}$, the generic fibre of $G_2$.
\end{theorem}

%This corrects the proof of \cite{Ant}, Lemma 3.3. 
This immediately implies that when $G$, $G_1$ and $G_2$ are all models of a same $K$-group scheme $G_{K}$ ($K$ being the fraction field of $R$) then the pushout of (\ref{eqIntro}) exists and is still a model of $G_{K}$ thus proving the existence of a lower bound for models of finite group schemes. This was already known in the commutative case (cf. \cite{RA}, Proposition 2.2.2). The same will be true for the quasi-finite case under the assumption that $G_{2,K}$ admits a finite and flat $R$-model:

\begin{theorem}(Cf. Theorem \ref{teoQFinito})
Assume $R$ is a complete discrete valuation ring and $G$, $G_1$, $G_2$ are quasi-finite and flat over $R$. If $\varphi_1$ is a model map and $G_{2,K}$ admits a finite and flat model then the pushout of (\ref{eqIntro}) in the category of affine $R$-group schemes exists. Moreover it is quasi-finite and flat and its generic fibre is isomorphic to $G_{2,K}$.
\end{theorem}

Using the fact that  $G_{2,K}$ always admits, when it is \'etale, a finite and flat model up to a finite extension of scalars we finally prove the following

\begin{corollary}(Cf. Corollary \ref{corQFinito2})
Assume $R$ is a complete discrete valuation ring and $G$, $G_1$, $G_2$ are quasi-finite and flat over $R$. Then if $\varphi_1$ is a model map  and $G_{2,K}$ is \'etale then  the pushout of (\ref{eqIntro}) in the category of  affine $R$-group schemes exists. Again it is quasi-finite and flat and its generic fibre is isomorphic to $G_{2,K}$.
\end{corollary} 

All the proofs rest on the computation of the pushout in the category of $R$-Hopf algebras. With the same techniques we briefly study in section \S \ref{sez:cok} the existence of cokernels in the category of affine $R$-group schemes where $R$ is any Dedekind ring. This will lead to a new and short proof of the the following:

\begin{corollary}(Cf. Corollary \ref{corEsisteCok})
Let $R$ be a Dedekind ring, $G$ and $H$ two finite and flat $R$-group schemes with $H$ a closed and normal $R$-subgroup scheme of $G$. Then the quotient $G/H$ exists in the category of $R$-affine group schemes.
\end{corollary} 

This holds over any base scheme and  is in fact a consequence of a much bigger theorem (Cf. \cite{SGA} Th\'eor\`eme 7.1).

\subsection{Notations and conventions} Every ring $A$ will be supposed to be associative and unitary, i.e. provided with a unity element denoted by $1_A$, or simply $1$ if no confusion can arise. However, unless stated otherwise, a ring will not be supposed to be commutative.  Every Dedekind ring, instead, will  always be supposed to be commutative. For an $R$-algebra $A$ the morphisms $u_A:R\to A$ and $m_A:A\otimes_R A\to A$ will always denote the unity and the multiplication morphisms (respectively). If moreover $A$ has an $R$-coalgebra structure then $\Delta_A:A\to A\otimes_R A$, $\varepsilon_A:A\to R$ will denote the comultiplication and the counity  respectively. Furthermore if $A$ has a $R$-Hopf algebra structure then $S_A:A\to A$ will denote  the coinverse. All the coalgebra structures will be supposed to be coassociative. Morphisms of $R$-algebras (resp. $R$-coalgebras, $R$-Hopf algebras) are $R$-module morphisms preserving $R$-algebra (resp. $R$-coalgebra, $R$-Hopf algebra) structure. We denote by $R$-$\mathcal{H}opf$ the category of associative and coassociative $R$-Hopf algebras while $R$-$\mathcal{H}opf_{ff}$ will denote the category of associative and coassociative $R$-Hopf algebras which are finite and flat as $R$-modules. When $R\to T$ is a morphism of commutative algebras, $M$ is a $R$-module, $X$ is a $R$-scheme, $f:M\to N$ is a $R$-module morphism and $\varphi:X\to Y$ a morphism of $R$-schemes then we denote by $M_T$, $X_T$, $f_T:M_T\to N_T$ and $\varphi_T:X_T\to Y_T$   respectively the $T$-module $M\otimes_RT$,  the $T$-scheme $X\times_{Spec(R)}Spec(T)$, the $T$-module morphism induced by $f$ and the $T$-morphism of schemes induced by $\varphi$.  When  $R$ is a Dedekind ring and $K$ its field of fractions then a $R$-morphism of schemes $\varphi:X\to Y$ is called a model map if generically it is an isomorphism, i.e. $\varphi_K:X_K\to Y_K$ is an isomorphism.

\section{Pushout of Hopf algebras}
%\subsection{The main construction}
\label{sez:poha}
In this section we first study the pushout of algebras over a commutative ring  $R$ then we  discuss the existence of the pushout in the category of  $R$-$\mathcal{H}opf_{ff}$ when $R$ is a complete discrete valuation ring. Let us first recall that  in a category $\mathcal{C}$ the pushout (see for instance \cite{ML}, III, \S 3) of a diagram
\begin{equation}\label{eqPush}\xymatrix{  & A\ar[ld]_{f}\ar[rd]^{g} & \\B & & C}\end{equation}
(where clearly $A,B,C$ are objects of $\mathcal{C}$ and $f,g$ morphisms in the same category) is an object of $\mathcal{C}$ that we denote $B\sqcup_A C$ provided with two morphisms $u:B\to B\sqcup_A C$, $v:C\to B\sqcup_A C$ such that $uf=vg$ and satisfying the following universal property: for any object $P$ of $\mathcal{C}$ and any two morphisms $u':B\to P$, $v':C\to P$ in $\mathcal{C}$ such that $u'f=v'g$ then there exists a unique morphism $p:B\sqcup_A C\to P$ making the following diagram commute:
$$
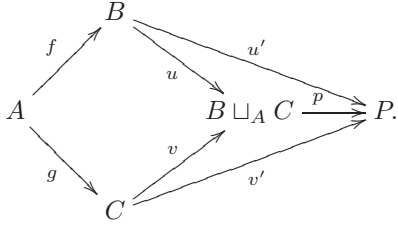
\xymatrix{  & B\ar[rd]_u \ar[rrd]^{u'} & & \\  A\ar[ur]^f \ar[dr]_g & & B\sqcup_A C\ar[r]^p & P.\\  & C\ar[ru]^v \ar[rru]_{v'} & & }$$
When $A$ is an initial object (provided it exists) of $\mathcal{C}$ then $B\sqcup_A C$ is the coproduct\footnote{The coproduct can be defined, however,  without assuming the existence of an initial object.} of $B$ and $C$ in $\mathcal{C}$.
%Let us recall the following lemma which follows directly from the universal property of the pushout:
%\begin{lemma}\label{lemPrimoPO}
%Let $A, B, C, D$ be objects of a category $\mathcal{C}$ provided with morphisms $A\to B$, $A\to C$, $A\to D$,  then 
%\begin{enumerate}
%\item $B\sqcup_AC\simeq C\sqcup_AB$,
%\item $B\sqcup_AA\simeq B$ (where $A\to A$ is the identity map),
%\item $(B\sqcup_AC)\sqcup_AD \simeq B\sqcup_A (C\sqcup_AD)$.
%\end{enumerate}
%\end{lemma}
When $\mathcal{C}$ is the category of commutative $R$-algebras then the pushout is given by the tensor product $B\otimes_A C$. This is not true anymore if $\mathcal{C}$ is the category of $R$-algebras,  (cf. Example \ref{exaCopro3} or create easier examples). However we can always find a pushout even when $\mathcal{C}$ is the category of $R$-algebras and it will be denoted by  $B\ast_A C$.  Before introducing, however, the pushout for non (necessarily) commutative $R$-algebras we recall the behavior of the tensor product over $R$.  We put ourselves in the following  situation:

\begin{notation}\label{notSiComincia0}
By  $R$ we will denote a commutative ring while  $A$, $B$ and $C$ will be $R$-algebras and $f:A\to B$, $g:A\to C$ two $R$-algebra morphisms. We also denote by $\rho_B:B\to B\otimes_R C$ and $\rho_C:C\to B\otimes_R C$ the morphisms sending respectively $b\mapsto b\otimes 1_C$ and $c\mapsto 1_B\otimes c$.
\end{notation}

\begin{proposition}\label{propUnivTensor}
Let $D$ be any $R$-algebra and $u:B\to D$, $v:C\to D$ two $R$-algebra morphisms such that $u\circ u_B=v \circ u_C$ and such that $u(b)v(c)=v(c)u(b)$ for all $b\in B, c\in C$. Then there exists a unique $R$-algebra morphism $t:B\otimes_R C\to D$ making the following diagram commute:
$$
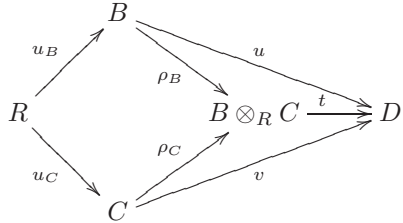
\xymatrix{ & & B\ar[rd]_{\rho_B}\ar[rrd]^{u} & \\ & R\ar[ru]^{u_B}\ar[rd]_{u_C} &  & B\otimes_R C\ar[r]^{t} & D \\ & & C\ar[ru]^{\rho_C}\ar[rru]_{v}  &}$$
\end{proposition}
\proof Cf. for instance \cite{FD}, I, \S 3 Proposition 3.2.
%BLAB The existence of the $R$-module morphism $t:B\otimes_AC \to D; b\otimes_A c\mapsto u(b)v(c)$ in ensured by the  $A$-balanced map $B\times C\to D; (b,c)\mapsto u(b)v(c)$. That $t$ is an $R$-algebra morphism follows from the property $u(b)v(c)=v(c)u(b)$ for all $b\in B, c\in C$.
\endproof

Unfortunately $B\otimes_AC$ behaves badly in general and one can observe that even $A\otimes_AA\simeq A$, as an $R$-algebra, is not a natural quotient of $A\otimes_R A$. So, instead, let us consider the following construction:

\begin{definition}\label{defQuozTensor}
We denote by $B\star_A C$, and we call it the star product of $B$ and $C$ over $R$, the quotient of $B\otimes_R C$ by the two-sided ideal generated by $A$, i.e. the ideal of $B\otimes_R C$ generated by the set $\{\rho_Bf(a)-\rho_Cg(a)\}_{a\in A}$. 
\end{definition}
It is an easy consequence the following universal property of the star product:
\begin{proposition}\label{propUnivStar}
Let $D$ be any $R$-algebra and $u:B\to D$, $v:C\to D$ two $R$-algebra morphisms such that $uf=vg$ and such that $u(b)v(c)=v(c)u(b)$ for all $b\in B, c\in C$. Then there exists a unique $R$-algebra morphism $t:B\star_A C\to D$ making the following diagram commute:
$$
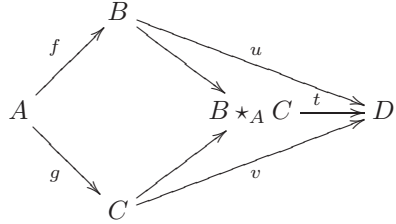
\xymatrix{ & & B\ar[rd]\ar[rrd]^{u} & \\ & A\ar[ru]^{f}\ar[rd]_{g} &  & B\star_A C\ar[r]^{t} & D \\ & & C\ar[ru]\ar[rru]_{v}  &}$$
\end{proposition}
\proof It is sufficient to take the $R$-algebra morphism $B\otimes_RC\to D$ and observe that it passes to the quotient.
\endproof

 The star product will only be used in Example \ref{exaCopro2} and  \S \ref{sez:cok}, so finally let us recall the  construction of the pushout of $R$-algebras: we follow essentially \cite{CPM}  1.7 and 5.1  with very few modifications in the exposition in order to obtain an easier to handle description. We  describe $A$, $B$ and $C$ giving their presentation as  $R$-algebras thus getting  $R\langle X_0;S_0\rangle$, $R\langle X_1;S_1\rangle$ and $R\langle X_2;S_2\rangle$ respectively, where $X_i$ is a generating set with relations $S_i$ ($i=0, 1, 2$). We recall that $R\langle X;S\rangle$ is to be intended as the $R$-algebra whose elements are all $R$-linear combinations of words on the set $X$ quotiented by the two-sided ideal generated by the relations in $S$.  Observe that  for  $y,z\in X$ we are not assuming $zy=yz$; if it is the case the information will appear in $S$. However for any $x\in X$ and any $r\in R$ we do assume $xr=rx$. For example the commutative $R$-algebra $R[x,y]/f(x,y)$ can be presented as $R\langle x,y;f(x,y)=0,xy=yx \rangle$. First we observe that the coproduct of $B$ and $C$ (i.e. the pushout of $B$ and $C$ over the initial object $R$) is given by the $R$-algebra $B\ast_R C:=R\langle X_1\cup X_2; S_1\cup S_2\rangle$ where the union is of course disjoint. Let us denote by $u:B\to B\ast_R C$ and $v:C \to B\ast_R C$ the canonical inclusions. Then the pushout of $B$ and $C$ over $A$ is given by the $R$-algebra 

\begin{equation}\label{eqAST}
B\ast_A C:=R\langle X_1\cup X_2; S_1\cup S_2\cup S_3\rangle
\end{equation} where $S_3$ consists on the relations given by  $uf(x)=vg(x)$ for every $x\in X_0$. Now we relate the pushout just described to the tensor product:

%\begin{lemma}\label{lemSuriet2}
%The universal morphism $\varphi: B\ast_RC\to B\otimes_RC$ is surjective.
%\end{lemma}
%\begin{proof} If it is not we can factor $\varphi$ as follows 
%$$\xymatrix{B\ast_RC\ar@{->>}[r]^{\varphi_2} & Q\ar@{^{(}->}[r]^{\varphi_1} & B\otimes_RC}$$ then we have the following commutative diagram
%$$\xymatrix{
%& B\ar[rd]_{u'}\ar[rrd]_{u''}\ar[rrrd]^{u} & & & \\R\ar[ru]\ar[rd] & & B\ast_RC\ar@{->>}[r]^{\varphi_2} & Q\ar@{^{(}->}[r]^{\varphi_1} & B\otimes_RC\\ & C \ar[ru]^{v'}\ar[rru]^{v''} \ar[rrru]_{v}& & & }$$
%We know that for any $b\in B$ and any $c\in C$ we have $u(b)v(c)=v(c)u(b)$ but this implies that $\varphi_1(u''(b)v''(c))=\varphi_1(v''(c)u''(b))$, but $\varphi_1$ is injective so finally  $u''(b)v''(c)=v''(c)u''(b)$, then we conclude by means of the universal property of the tensor product (cf. Proposition \ref{propUnivTensor}). 
%\end{proof}

\begin{lemma}\label{lemSuriet}
Assume that $B=R\langle X_1; S_1\rangle$ and $C=R\langle X_2; S_2\rangle$. Then $B\otimes_RC$ can be presented as $R\langle X_1\cup X_2; S_1\cup S_2, \{zy=yz\}_{z\in X_1,y\in X_2}\rangle$
 thus becoming a quotient of $R\langle X_1\cup X_2; S_1\cup S_2\rangle=B\ast_RC$.
\end{lemma}
\proof Let $D$ be an $R$-algebra provided with $R$-algebra morphisms $m:B\to D$ and $n:C\to D$ such that $m\circ u_B=n\circ u_C$, and assume moreover that $m(b)n(c)=n(c)m(b)$ for all $b\in B, c\in C$. Let us denote by $u:B\to B\ast_RC$ and  $v:C\to B\ast_RC$ the canonical morphisms and by $\lambda:B\ast_R C\to D$ the universal morphism making the following diagram commute:
$$\xymatrix{ &B\ar[dr]_u\ar[rrd]^m & & \\R\ar[ur]^{u_B}\ar[dr]_{u_C} & &B\ast_RC\ar[r]^{\lambda} & D \\ &C\ar[ur]^v\ar[rru]_n & & }$$
By assumption $\lambda u(z)\lambda v(y)=\lambda v(y)\lambda u(z)$ so $ u(z)v(y)- v(y) u(z)\in ker(\lambda)$ hence $\lambda$ factors through $R\langle X_1\cup X_2; S_1\cup S_2, \{zy=yz\}_{z\in X_1,y\in X_2}\rangle$ providing it with the universal property stated in Proposition \ref{propUnivTensor} and this is enough to conclude.
%have a natural morphism $\delta: B\otimes_RC\to R\langle X_1\cup X_2; S_1\cup S_2, \{zy=yz\}_{z\in X_1,y\in X_2}\rangle$ which is surjective since the composition $$B\ast_R  C \to B\otimes_RC\to R\langle X_1\cup X_2; S_1\cup S_2, \{zy=yz\}_{z\in X_1,y\in X_2}\rangle$$ is surjective. It is easily seen to be injective working on generators.
\endproof

Let $q:R\to T$ be a $R$-commutative algebra. When $f=f(x_1, .., x_n)\in R[x_1, .., x_n]$ we denote by $q_*(f)$ the polynomial in $T[x_1, .., x_n]$ whose coefficients are the image in $T$ by $q$ of the coefficients of $f$, i.e. the image of $f$ through the morphism $q_{\ast}:R[x_1, .., x_n]\to T[x_1, .., x_n]=R[x_1, .., x_n]\otimes_RT$. Now take $R\langle X, S\rangle $: by an abuse of notation we denote by $q_*(S)$ the set of relations $\{q_*(s_i)=0\}$ on the set $X$. In Lemma \ref{lemBaseChange} we  observe that the pushout is stable under base change.

\begin{lemma}\label{lemBaseChange}
Let $q:R\to T$ be a $R$-commutative algebra and $R\langle X;S\rangle$ any $R$-algebra, then 
\begin{enumerate}
\item $R\langle X;S\rangle\otimes_RT\simeq T\langle X;q_*(S)\rangle$,
\item $(B\ast_AC)\otimes_RT \simeq (B\otimes_RT)\ast_{(A\otimes_RT)} (C\otimes_RT)$.
\end{enumerate}
\end{lemma}
\proof
As a commutative $R$-algebra, $T$ is isomorphic to $R[\{y_i\}]/(\{f_r\})$ where $\{y_i\}$ is a set of generators and $\{f_r\}$ a set of polynomials in the variables $\{y_i\}$ with coefficients in $R$. So by Lemma \ref{lemSuriet} $R\langle X;S\rangle\otimes_RT$ is isomorphic to $R\langle X\cup \{y_i\}; S\cup \{f_r=0\} \cup \{y_iy_j=y_jy_i\}\cup \{xy_i=y_ix\}_{x\in X} \rangle$ which is isomorphic to $R\langle X\cup \{y_i\}; q_*(S)\cup \{f_r=0\} \cup \{y_iy_j=y_jy_i\}\cup \{xy_i=y_ix\}_{x\in X} \rangle$ and the latter is isomorphic to $T\langle X;q_*(S)\rangle$ since $T$ commutes with $X$ and this proves 1. Let us describe $A$, $B$ and $C$ as   $R\langle X_0;S_0\rangle$, $R\langle X_1;S_1\rangle$ and $R\langle X_2;S_2\rangle$ respectively. As a consequence of point 1 we have  $A\otimes_RT\simeq T\langle X_0;q_*(S_0)\rangle$,  $B\otimes_RT\simeq T\langle X_1;q_*(S_1)\rangle$, $C\otimes_RT\simeq T\langle X_2;q_*(S_2)\rangle$ and $(B\ast_AC)\otimes_RT\simeq T\langle X_1\cup X_2; q_*(S_1\cup S_2\cup S_3)\rangle$ where $S_3$ is as described in (\ref{eqAST}).   But $(B\otimes_RT)\ast_{(A\otimes_RT)} (C\otimes_RT)$ is also isomorphic to the latter which enables us to conclude.
\endproof

\begin{notation}\label{notFM}
When $R$ is a  Dedekind ring and $M$ an $R$-module, let us denote by $q:M\twoheadrightarrow F(M)$ the unique quotient (cf. \cite{EGAIV2}, Lemme (2.8.1.1)) of $M$ which is $R$-flat and such the induced map $q_K:M_K\to F(M)_K$ is an isomorphism.
\end{notation} 

Let us analyze a few examples whose importance will be clear in the following sections:

%\begin{example}\label{exaCopro}
%Let us set $A:=R[x]/f(x)$, thus a commutative $R$-algebra, for some polynomial $f(x)\in R[x]$. It is known that $A\otimes_RA=\frac{R[x,y]}{f(x),f(y)}$ and it is the pushout of $A$ and $A$ over $R$ (hence the coproduct) in the category of commutative $R$-algebras. Now we  describe $A\ast_RA$: it is presented by $R\langle x,y;f(x)=f(y)=0  \rangle$ thus $xy\neq yx$ then the natural surjection (cf. Lemma \ref{lemSuriet}) $A\ast_RA\twoheadrightarrow A\otimes_RA$ is not an isomorphism. Observe also that $A\ast_RA$ is not even finite as an $R$-module as $x, y, xy, xyx, xyxy, xyxyx, ...$ are all $R$-linearly independent.
%\end{example}

\begin{example}\label{exaCopro2} Let $f:A\to B$ and $g:A\to R$ be morphisms of $R$-algebras then the canonical morphism $\varphi:B\ast_AR \to B\star_A R$ is an isomorphism. Indeed  we observe that $B\ast_RR=B$  and that the canonical morphisms  $u:B\to B\ast_RR$ and $v:R\to B\ast_RR$ are nothing else but $Id_B$ and the  unit morphism $u_B$ respectively, then for any $b\in B$ and any $r\in R$ we have $u(b)v(r)=v(r)u(b)$. Hence denoting by  $u':B\to B\ast_AR$ and $v':R\to B\ast_AR$  the canonical morphisms we also have  $u'(b)v'(r)=v'(r)u'(b)$ as $u'=\lambda u$ and $v'=\lambda v$ where $\lambda:B\ast_RR \to B\ast_AR$ is the universal morphism. Then $\varphi$ can be inverted according to Proposition \ref{propUnivStar}. Observe that $B\ast_AR$ is finite as an $R$-module if $B$ is finite (it is indeed a quotient of $B$).
\end{example}

\begin{example}\label{exaCopro3} Let $R$ be a discrete valuation ring with uniformising element $\pi$. Let us fix a positive integer $p$ and let us set $A:=R[x]/x^p$, $B:=R[y]/y^p$ and $C:=R[z]/z^p$ (thus commutative $R$-algebras). Consider the morphisms $f:A\to B, x\mapsto \pi^ny$ and $g:A\to C, x\mapsto \pi^mz$ where $m> n> 0$ are integers. Then $B\ast_AC=R\langle y,z;y^p=z^p=0, \pi^mz=\pi^ny  \rangle$. Observe that, as an $R$-module, $B\ast_AC$ is not  flat as $\pi^{n}(\pi^{m-n} z-y)=0$ thus $\pi^{m-n} z-y$ is a $R$-torsion element. However if we add the relation $\pi^{m-n} z=y$ then we eliminate torsion from $B\ast_AC$ and what we obtain is (cf. Notation \ref{notFM}) $F(B\ast_AC)=R\langle y,z;y^p=z^p=0, \pi^{m-n}z=y  \rangle=R[z]/z^p$ thus finitely generated  and flat and, in this particular case, it is isomorphic to $F(B\otimes_AC)$.
%$xy=0=yx$ to $B\ast_AC$ (i.e. we take a suitable quotient of  $B\ast_AC$) we obtain $R\langle y,z;y^p=z^p=0, \pi^mz=\pi^ny,zy=0=yz  \rangle$

%st relation can be rewritten as $a^{m-n}z=y$  hence $B\ast_AC=R\langle z;z^p=0  \rangle=C$ and the canonical morphisms to $B\ast_AC $ are thus given by $u':B\to B\ast_AC, y\mapsto a^{m-n}z$ and $v':C\to B\ast_AC, z\mapsto z$ which is nothing but the identity on $C$. Again we have $B\ast_AC\simeq B\otimes_AC$.
\end{example}

%When $R$ is a discrete valuation ring of uniformising paramter $\pi$ then the field of fractions $K:=Frac(R)$ can be expressed as $K=R[x]/(\pi x-1)=R\langle x; \pi x=1\rangle$ and $R\langle X; S\rangle\otimes_RK=R\langle X, x; S, \pi x=1, \{xz=zx\}_{z\in X} \rangle=K\langle X; S\rangle$. Hence it is easy to check that $(B\ast_AC)\otimes_R K\simeq (B\otimes_RK)\ast_{(A\otimes_RK)}(C\otimes_RK)$. More generally we have the following
%\begin{lemma}\label{lemEstensScalar}
%Let $R$ be an integral domain and $K$ its field of fractions. Let $A,B, C$ be $R$-algebras as in the following commutative diagram 

%\begin{equation}\label{eqZero}\xymatrix{ & & B\\ R\ar[r]^{u_A}\ar[rru]^{u_B}\ar[rrd]_{u_C} & A\ar[ur]_f \ar[dr]^g & \\ & & C}\end{equation} 

%Then $(B\ast_AC)\otimes_R K\simeq (B\otimes_RK)\ast_{(A\otimes_RK)}(C\otimes_RK)$.
%\end{lemma}
%\proof As in the construction preceding the statement it is sufficient to observe that $R\langle X;S\rangle\otimes_RK=K\langle X;S\rangle$.
%\endproof
%{\color{red}Questo corollario e` da cambiare in quanto falsissimo}

%It worths recalling the following

%\begin{lemma}\label{lemStamani}
%Let $R$ be a Dedekind ring with fraction field $K$ and for any prime $p\in R$ let us set $k(p):=R/pR$; let $M$ be a torsion free $R$-module of finite rank . Then for every prime $p\in R$ we have $dim_{k(p)}(M\otimes_Rk(p))\leq dim_K(M\otimes_RK)$.
%\end{lemma}
%\proof Well known and easy to check.
%\endproof

The following well known result will be used several times in this paper:

\begin{theorem}\label{teoKapl}
Let $R$ be a complete discrete valuation ring with fraction field $K$ and residue field $k$. Let $M$ be a torsion-free $R$-module of finite rank $r$ (i.e. $r:=dim_K(M\otimes_RK)<+\infty$). Then $M\simeq_{R\text{-}mod} K^{\oplus r-s}\oplus R^{\oplus s}$, where  $s=dim_k(M\otimes_Rk)$. 
\end{theorem}
\proof
This is \cite{KAP}, Chapter 16, Corollary 2, 
\endproof

Theorem \ref{teoKapl} is not true when $R$ is not complete (cf. \cite{KAP}, Theorem 19) and this is why we will often need to restrict to complete discrete valuation rings. The following lemma is crucial in this paper:

\begin{lemma}\label{lemEstensScalar}
Let $R$ be a complete discrete valuation ring and assume that $f:A\to B$ and $g:A\to C$ are $R$-algebra morphisms where furthermore $g_K:A_K\to C_K$ is  an isomorphism. Then if $A$, $B$ and $C$ are finitely generated and flat as $R$-modules then the same holds for $F(B\ast_AC)$. Moreover the canonical $R$-algebra morphism $B\to F(B\ast_AC)$ induces an isomorphism $B_K\to F(B\ast_AC)_K$.
\end{lemma}
\proof Let $\pi$ be an uniformising element of $R$ and let $K$ and $k$ be the fraction and residue fields respectively. As usual let us present by $R\langle X_0;S_0\rangle$, $R\langle X_1;S_1\rangle$ and $R\langle X_2;S_2\rangle$ respectively the $R$-algebras $A$, $B$ and $C$ where for $X_0, X_1, X_2$  we take respectively  bases of $A$, $B$, $C$ as $R$-modules minus the identity elements so that the cardinality  $x_0, x_1, x_2$ of those sets is the rank of $A, B$ and $C$ (respectively) minus one; of course $x_0=x_2$. 
%Again by \cite{EGAIV2}, Lemme 2.8.1.1 we observe that $F(B\ast_AC)$ is isomorphic to the unique quotient of $B\ast_RC$ which is $R$-flat and generically isomorphic to $B\otimes_RK$: indeed $B\otimes_RK \simeq (B\otimes_RK)\ast_{(A\otimes_RK)} (C\otimes_RK)$ is isomorphic to $(B\ast_AC)\otimes_RK$ by Lemma \ref{lemBaseChange}. 
Now $B\ast_AC=R\langle X_1\cup X_2;S_1\cup S_2\cup S_3\rangle$ where $S_3$ is as described in (\ref{eqAST}). In particular $S_3$ is a set made of $x_0$ $R$-linear relations relating the $x_2$ elements of $X_2$ and the $x_1$ elements of $X_1$. As in Example \ref{exaCopro3} the information of being $R$-torsion (if any) is contained in the set $S_3$, so if we want to cut out $R$-torsion we need to add another set of $x_0$ relations $S_4$ obtained as follows: for each relation $(s=0)\in S_3$ add the relation $(t=0)$ to $S_4$ where $s=\pi^vt$ and $t$ has at least one coefficient equal to an invertible element of $R$. Thus $F(B\ast_AC)= R\langle X_1\cup X_2;S_1\cup S_2\cup S_3\cup S_4\rangle$. But since relations in $S_3$ are automatically satisfied if we add $S_4$ then $F(B\ast_AC)= R\langle X_1\cup X_2;S_1\cup S_2\cup S_4\rangle$. 
%Since $F(B\ast_AC)$ is torsion-free and of finite rank then $F(B\ast_AC)\otimes_Rk$ is a finite $k$-module of dimension  
%\begin{equation}\label{eqVoltaBuona}
%dim_k(F(B\ast_AC)\otimes_Rk)\leq r_B.
%\end{equation} 
Now Lemma \ref{lemBaseChange}, point 1, implies that $F(B\ast_AC)\otimes_Rk=k\langle X_1\cup X_2;q_*(S_1\cup S_2\cup S_4)\rangle$ where $q:R\twoheadrightarrow k$ is the canonical surjection so $F(B\ast_AC)\otimes_Rk$ is the quotient of  $k\langle X_1\cup X_2;q_*(S_1\cup S_2)\rangle$ by the two-sided ideal generated by the relations $q_*(S_4)$. But in $k\langle X_1\cup X_2;q_*(S_1\cup S_2)\rangle$ the elements of the set $X_1\cup X_2$ are $x_1+x_2$ $k$-linearly independent vectors then if we add the $x_0=x_2$ $k$-linear relations $q_*(S_4)$ what remains is a set of at least $x_1=rk(B)-1$ $k$-linearly independent elements which become $rk(B)$ if we add $1_B$. Combining this with %(\ref{eqVoltaBuona}) we obtain $dim_k(F(B\ast_AC)\otimes_Rk)= r_B$ and 
Theorem \ref{teoKapl} we obtain that $F(B\ast_AC)$ is a finitely generated $R$-free module, as required, as $dim_k(F(B\ast_AC))\otimes_Rk=dim_K(F(B\ast_AC))\otimes_RK$. The last assertion follows easily from Lemma \ref{lemBaseChange}, point 2.

%By Lemma \ref{lemBaseChange} we know that $\varphi$ is generically an isomorphism. Then  $\varphi$ induces a morphism $F(\varphi):F(B\ast_AC)\to F(B\otimes_AC)$  by \cite{EGAIV2}, (2.8.3) which is surjective and generically isomorphic as $\varphi$ is  (cf. Lemma \ref{lemSuriet}).  Hence we conclude, again,  by \cite{EGAIV2}, Lemme (2.8.1.1).
\endproof

\begin{remark}\label{lemNonDipende}
The construction in Lemma \ref{lemEstensScalar} does not depend on $A$. That means that if we take $A'$, $f':A'\to B$ and $g':A'\to C$ satisfying similar assumptions then $F(B\ast_AC)\simeq F(B\ast_{A'}C)$. Indeed, again by \cite{EGAIV2}, Lemme 2.8.1.1, we observe that $F(B\ast_AC)$ is isomorphic to the unique quotient of $B\ast_RC$ which is $R$-flat and whose tensor over $K$ gives $B_K$; but the same property is satisfied  by  $F(B\ast_{A'}C)$ hence we conclude by unicity.
\end{remark}

%Let notation be as in \ref{notSiComincia} 

% In \cite{BB} (6.3) Barr and Beck construct the coproduct in the category of supplemented $R$-algebras while in \cite{JWB}, 4.1 Jordan constructs, in the same category, the pushout in the category of supplemented algebras whenever  $A\to B$ and $A\to C$ are monomorphisms. Let us briefly recall both constructions:

%Let us show that whenever the objects involved are $R$-Hopf algebras and not just supplemented algebras then Jordan's construction allows us to obtain the pushout in the category of $R$-Hopf algebras.
%{\color{red} occhio qui devi cambiare piu` o meno tutto}
%Let $A,B, C$ be $R$-Hopf?? algebras as in the following commutative diagram 

%\begin{equation}\label{eqZero}\xymatrix{ & & B\\ R\ar[r]^{u_A}\ar[rru]^{u_B}\ar[rrd]_{u_C} & A\ar[ur]_f \ar[dr]^g & \\ & & C}\end{equation} 

\begin{proposition}\label{propPOHA}
Let $R$ be any commutative ring. Then the pushout in the category of  $R$-bialgebras exists.
\end{proposition}
\proof We follow\footnote{In \cite{LJM}, however, Lemaire uses different notations.}  \cite{LJM}, Chapitre 5, \S 5.1, Proposition. Consider the diagram 
\begin{equation}\label{eqZero}\xymatrix{ & & B\\ R\ar[r]^{u_A}\ar[rru]^{u_B}\ar[rrd]_{u_C} & A\ar[ur]_f \ar[dr]^g & \\ & & C}\end{equation} 
where we assume that $A$, $B$ and $C$ are $R$-bialgebras and the arrows are $R$-Hopf algebra morphisms. Let $D:=B\ast_A C$  be the pushout of the diagram in the category of $R$-algebras and let $m_{D}$ and $u_{D}$ be respectively the multiplication and the unit morphism. Then we need to provide $D$ with a comultiplication $\Delta_{D}$ and a counit $\varepsilon_{D}$  such that $({D},m_{D},u_{D},\Delta_{D},\varepsilon_{D})$ is a $R$-bialgebra. We describe how to construct $\Delta_{D}$, the costruction of $\varepsilon_D$ being easier. The rest will be standard verification over complicated diagrams. The existence of $\Delta_{D}$ is explained in the following diagram, taking into account the universal property of $D$:
$$\xymatrix{& B\ar[rd]^u\ar[rrr]^{\Delta_B} & & & B\otimes_R B \ar[rd]^{u\otimes u} & \\A\ar[ru]^f\ar[rd]_g\ar@/^/[rrr]^(.3){\Delta_A} & & D\ar@{-->}@/_/[rrr]_(.7){\Delta_D} & A\otimes_R A \ar[ru]^{f\otimes f}\ar[rd]_{g\otimes g}& & D\otimes_R D \\ & C\ar[ru]_v\ar[rrr]_{\Delta_C} & & & C\otimes_R C \ar[ru]_{v\otimes v} & } $$
\endproof

\begin{remark}\label{remANTIPODO}
Notation being as in Proposition \ref{propPOHA} one observes that we can define a $R$-module morphism $S_D:D\to D$, candidate to be a coinverse,  as follows: first construct the opposite algebras $A^{op}$, $B^{op}$, $C^{op}$, $D^{op}$, the opposite morphisms $f^{op}$, $g^{op}$, $u^{op}$, $v^{op}$ and the morphisms of $R$-algebras $S'_A$, $S'_B$, $S'_C$, induced by the $R$-algebra anti-morphisms $S_A$, $S_B$, $S_C$. Then the existence of $S'_D$ follows from the following diagram
$$
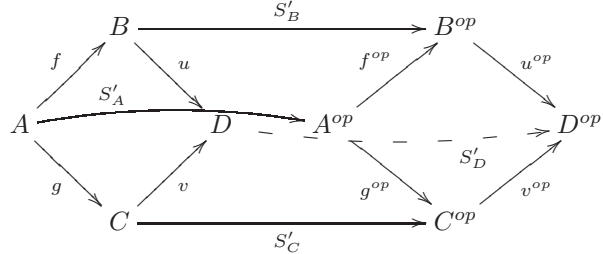
\xymatrix{& B\ar[rd]^u\ar[rrr]^{S'_B} & & & B^{op} \ar[rd]^{u^{op}} & \\A\ar[ru]^f\ar[rd]_g\ar@/^/[rrr]^(.3){S'_A} & & D\ar@{-->}@/_/[rrr]_(.7){S'_D} & A^{op} \ar[ru]^{f^{op}}\ar[rd]_{g^{op}}& & D^{op} \\ & C\ar[ru]_v\ar[rrr]_{S'_C} & & & C^{op} \ar[ru]_{v^{op}} & } $$
exploiting the universal property of $D$ then $S_D$ is the anti-morphism induced by $S'_D$. However $S_D$ may fail to be a coinverse for $D$ as $m_D\circ(S_D\otimes id_D)\circ m_D$ may not be equal to $u_D\circ \varepsilon_D$ (same for $m_D\circ(id_D \otimes S_D)\circ m_D$).
\end{remark}

In order to have an explicit description for $S_D$, constructed in Remark \ref{remANTIPODO}, set, as usual, $A=R\langle X_0;S_0\rangle$, $A=R\langle X_1;S_2\rangle$ and $C=R\langle X_2;S_2\rangle$ so  $D=R\langle X_1\cup X_2;S_1\cup S_2\cup S_3\rangle$ where $S_3$ is as described in (\ref{eqAST}); it is sufficient to set $S_D(x_1):=S_B(x_1)$ for any $x_1\in X_1$, $S_D(x_2):=S_C(x_2)$ for any $x_2\in X_2$  for any $x_1\in X_1$, $S_D(x_1x_2):=S_D(x_2)S_D(x_1)$ and $S_D(x_2x_1):=S_D(x_1)S_D(x_2)$. It is well defined and is by construction an anti-isomorphism for the $R$-algebra $D$. A similar construction gives  an explicit description of $\Delta_D$, taking into account that $\Delta_D$ is a morphism of $R$-algebras and not an anti-morphism.

\begin{corollary}\label{corPOHA}
Let $R$ be a complete discrete valuation ring and assume that $f:A\to B$ and $g:A\to C$ are $R$-algebra morphisms where furthermore $g_K:A_K\to C_K$ is  an isomorphism. Then $F(B\ast_AC)$ has a natural structure of $R$-Hopf algebra.  If moreover $A$, $B$ and $C$ are finitely generated and flat as $R$-modules then $F(B\ast_AC)$ is the pushout of $B$ and $C$ over $A$ in $R$-$\mathcal{H}opf_{ff}$.
\end{corollary}
\proof By \cite{EGAIV2}, (2.8.3) and of course Proposition \ref{propPOHA} we obtain that $F(B\ast_AC)$ has a natural structure of $R$-bialgebra. We need to prove the existence of a coinverse $S_{F(B\ast_AC)}$ that gives $F(B\ast_AC)$ a natural structure of $R$-Hopf algebra. So let us take for $D:=B\ast_AC$ the $R$-module morphism $S_{D}$ defined in Remark \ref{remANTIPODO}. This morphism induces (by \cite{EGAIV2}, Lemme 2.8.3) an $R$-module morphism $S_{F(D)}:F(D)\to F(D)$ which is the required coinverse: indeed 
\begin{equation}\label{eqAntipo} 
m_{F(D)}\circ(S_{F(D)}\otimes id_{F(D)})\circ m_{F(D)}-u_{D}\circ \Delta_{D}
\end{equation}  
is the zero map $0_D$ and this is clear since $F(D)\subset B_K$ and (\ref{eqAntipo}) tensored over $K$ gives rise to the equality
$$m_{B_K}\circ(S_{B_K}\otimes id_{B_K})\circ m_{B_K}=u_{B_K}\circ \Delta_{B_K}$$
which holds as $B_K$ is a $K$-Hopf algebra. The same is true for  $m_{F(D)}\circ(id_{F(D)} \otimes S_{F(D)})\circ m_{F(D)}$.  Finally  $F(B\ast_AC)$ is finitely generated and flat as an $R$-module when $A$, $B$ and $C$ are: this is Lemma \ref{lemEstensScalar}.
%{\color{red}occhio!!} Remark \ref{remEsisteq}.
\endproof

\begin{remark}
Notation being as in Proposition \ref{propPOHA},  we observe that   $B\ast_AC$ is  cocommutative if $A$, $B$, $C$ are. The same conclusion holds, then, for $F(B\ast_AC)$ in Corollary \ref{corPOHA}. Moreover observe that if $A$, $B$, $C$ are commutative then $F(B\ast_AC)$ is commutative too since it is contained in $B_K$. So in particular in this case  $F(B\ast_AC)\simeq F(B\otimes_AC)$, as it happend in Example \ref{exaCopro3}.
\end{remark}

\section{Pushout of group schemes}
\label{sez:PushOut}
In this section $R$ is any complete discrete valuation ring with fraction and residue fields respectively denoted by $K$ and $k$.
\subsection{The finite case} \label{sez:FinitePO} Let  $M=Spec(B)$ and $N=Spec(C)$ be finite and flat $R$-group schemes, i.e. $B$ and $C$ are free over $R$  and finitely generated  as $R$-modules. Let us  assume that there is a $K$-group scheme morphism $\psi:M_{K}\to N_{K}$. An \emph{upper bound} for $M$ and $N$ is a finite and flat $R$-group scheme $U$, provided with a model map $U\to M$ and a $R$-grouop scheme morphism $\varphi:U\to N$ which generically coincides with $\psi:M_{K}\to N_{K}$. A \emph{lower bound} for $M$ and $N$ is a finite and flat $R$-group scheme $L$, provided with a model map $N\to L$ and a $R$-grouop scheme morphism $\delta:M\to L$ which generically coincides with $\psi:M_{K}\to N_{K}$.  The construction of un upper bound is easy: it is sufficient to set $U$ as the schematic closure of $M_{K}$ in $M\times N$ through the canonical closed immersion $M_{K}\hookrightarrow M_{K}\times N_{K}$ (and this holds when the base  is any Dedekind scheme). Now consider the following commutative diagram
\begin{equation}\label{eqUno}\xymatrix{ & M\\  U\ar[ur]^m \ar[dr]_n & \\  & N}\end{equation} 
where $U=Spec(A)$ is any upper bound. We are now going to study the existence of a pushout $M\sqcup_U N$ in the category of finite and flat $R$-group schemes. We  prove the following

\begin{lemma}\label{lemFinito}
The pushout of (\ref{eqUno}) in the category of finite and flat $R$-group schemes exists. Moreover  $M\sqcup_U N$ is a lower bound for $M$ and $N$.
\end{lemma}
\proof
Notation being as in the beginning of this section, we have the following diagram of commutative  $R$-Hopf algebras:
$$\xymatrix{ & A & \\  B\ar[ru]  &  & C\ar[lu] }$$
which, dualizing, gives rise to the following diagram of cocommutative (but possibly non commutative) $R$-Hopf algebras: 
$$\xymatrix{ & A^{\vee}\ar[ld]_{f}\ar[rd]^g & \\  B^{\vee}  &  & C^{\vee}. }$$
Let us consider the cocommutative $R$-Hopf algebra $F(B^{\vee}\ast_{A^{\vee}}C^{\vee})$ constructed in Corollary \ref{corPOHA}. Now we take the spectrum of its dual  $P:=Spec(F(B^{\vee}\ast_{A^{\vee}}C^{\vee})^{\vee})$. First of all we observe that $F(B^{\vee}\ast_{A^{\vee}}C^{\vee})^{\vee}$ is commutative as $F(B^{\vee}\ast_{A^{\vee}}C^{\vee})$ is cocommutative so that taking its spectrum does make sense. It remains to prove that the commutative diagram
$$\xymatrix{ & M\ar[rd] & \\  U\ar[ur]^m \ar[dr]_n & & P\\  & N\ar[ur] & }$$
is in fact a pushout in the category of finite and flat $R$-group schemes. But this follows from the fact that $F(B^{\vee}\ast_{A^{\vee}}C^{\vee})$ is a pushout in $R$-$\mathcal{H}opf_{ff}$. That $M\sqcup_U N:=P$ is a lower bound for $M$ and $N$ is also clear by construction.
\endproof

\begin{theorem}\label{teoFinito}
The pushout of (\ref{eqUno}) in the category of affine $R$-group schemes exists. 
\end{theorem}
\proof
Consider the commutative diagram
$$\xymatrix{ & M\ar[rd]\ar[rrd]^u & & \\  U\ar[ur]^m \ar[dr]_n & & P & Q\\  & N\ar[ur] \ar[urr]_v & & }$$
where $P$ is the pushout of (\ref{eqUno}) in the category of finite and flat $R$-group schemes constructed in Lemma \ref{lemFinito} and $Q=Spec(D)$ is any affine $R$-group scheme. We are going to show that $P$ is also the pushout of (\ref{eqUno}) in the category of affine $R$-group schemes. Let us factor $u$ through $M':=Spec(B')$ via the morphisms $u':M\to M'$ and $i:M'\to Q$ where $i$ is a closed immersion and $u'$ is a schematically dominant morphism (i.e. the induced morphism $B'\to B$ is injective) so that  $M'$ is a finite and flat $R$-group scheme since $M$ is. Likewise we factor $v$ through the finite and flat $R$-group scheme $N':=Spec(C')$ via the schematically dominant morphism $v':N\to N'$ and the closed immersion $j:N'\to Q$. Now consider the finite $R$-group scheme (it needs not be flat a priori) $M'\times_{Q}N'$ and the natural closed immersions $i':M'\times_{Q}N'\hookrightarrow M'$ and $j':M'\times_{Q}N'\hookrightarrow N'$. So let us denote by $r:U\to M'\times_{Q}N'$ the universal morphism, then we have the following commutative diagram
$$\xymatrix{ & M\ar[rr]^{u'} & & M'\ar[rd]^i &\\  U\ar[ur]^m \ar[dr]_n \ar[rr]^r& & M'\times_{Q}N' \ar[ur]^{i'}\ar[dr]_{j'}&  & Q\\  & N \ar[rr]_{v'} & & N'\ar[ur]_{j} &}$$
and in particular we have the following commutative diagram of $R$-algebras

$$\xymatrix{A & B'\ar@{^{(}->}[l]\ar@{->>}[d]\\ & B'\otimes_DC'\ar[ul]}$$ 
where $B'\to B'\otimes_DC'$ is surjective but also injective since $B'\hookrightarrow A$ is (recall that $m:U\to M$ is a model map). Hence $i'$ is an isomorphism so there exists a universal morphism $P\to N'$ by Lemma \ref{lemFinito} and we conclude. 
\endproof

\begin{remark}\label{remNonDipendeFinito} If $U'$ is any other upper bound for $M$ and $N$ then there is a canonical isomorphism $M\sqcup_U N\simeq M\sqcup_{U'} N$:  this is a consequence of  Remark \ref{lemNonDipende}. However this does not mean that the lower bound  for $M$ and $N$ is unique, which is clearly not true in general.
\end{remark}

\subsection{The quasi-finite case}
\label{sez:qf}
 
Let  $M=Spec(B)$ and $N=Spec(C)$  be quasi-finite (by this we will always mean  affine and of finite type over $R$, with finite special and generic fibers) and flat $R$-group schemes. It is known (see \cite{BLR}, \S 7.3 p.179) that any  quasi-finite $R$-group scheme  $H$  has a finite part $H_f$, that is a open and closed subscheme of $H$ which consists of the special fibre $H_k$ and of all points of the generic fibre which specialize to the special fibre. It is thus flat over $R$ if $H$ is. 
\begin{remark}
If $H=Spec (A)$ is a quasi-finite and flat $R$-group scheme then its finite part coincides with $Spec(A^{\vee\vee})$ where $A^{\vee\vee}$ is the double dual of $A$: this follows from the fact that $A\simeq_{R\text{-}mod} K^{\oplus t}\oplus R^{\oplus s}$ (cf. Theorem \ref{teoKapl}) so $A^{\vee}\simeq_{R\text{-}mod} R^{\oplus s}$ is an $R$-Hopf algebra and not just an $R$-algebra.  Hence the canonical surjection $A\twoheadrightarrow A^{\vee\vee}$ gives the desired closed immersion $H_f\hookrightarrow H$ of group schemes. However this fact will not be necessary in the remainder of this paper.
\end{remark} 

 Let us  assume that there is a $K$-group scheme morphism $\psi:M_{K}\to N_{K}$. We define upper and lower bounds exactly as in the finite case. One can easily construct an upper bound $U$ for $M$ and $N$ simply proceeding as in \S \ref{sez:FinitePO}. So $U$ will be in general a quasi-finite and flat $R$-group scheme. For the lower bound it will be a little bit more complicated. So consider again the commutative diagram
(\ref{eqUno})
where $U=Spec(A)$ is any upper bound. We are going to study the existence of a pushout $M\sqcup_U N$ in the category of affine $R$-group schemes. We  prove the following

\begin{theorem}\label{teoQFinito}
Assume that $N_{K}$ admits a finite and flat model over $R$. Then the pushout of (\ref{eqUno}) in the category of affine $R$-group schemes exists. Moreover  $M\sqcup_U N$ is a lower bound for $M$ and $N$.
\end{theorem}
\proof Let $N'$ denote a finite and flat $R$-model for $N_{K}$, i.e. a finite and flat $R$-group scheme whose generic fibre is isomorphic to $N_{K}$. Consider the finite part $M_f$ and $N_f$ of, respectively, $M$ and $N$. Compose the closed immersion $M_{f,K}\to M_{K}$ with $\psi_{K}:M_{K}\to N_{K}$ thus obtaining a morphism $M_{f,K}\to N_{K}$. By Theorem \ref{teoFinito} we construct a lower bound $L_1$ for $M_f$  and $N'$, which is finite and flat over $R$, generically isomorphic to $N_{K}$,  then it is already a lower bound for $M$ and $N'$. Considering the closed immersion $N_{f,K}\to N_{K}$ we also construct a lower bound $L_2$ for $N_f$  and $N'$, which is finite and flat over $R$, generically isomorphic to $N_{K}$, then it is already a lower bound for $N$ and $N'$. So a lower bound $L$ for $L_1$ and $L_2$ (which are generically isomorphic) exists by Theorem \ref{teoFinito} and is also a lower bound for $M$ and $N$. We still need to compute the pushout of $M$ and $N$ over $U$: let us set $U_f:=Spec(A_f)$, $M_f:=Spec(B_f)$, $N_f:=Spec(C_f)$ and $L:=Spec(D)$. Consider the natural $R$-bialgebra morphism (cf. Proposition \ref{propPOHA}) $B_f^{\vee}\ast_{A_f^{\vee}}C_f^{\vee}\to D^{\vee}$ and factor it as follows
$$\xymatrix{B_f^{\vee}\ast_{A_f^{\vee}}C_f^{\vee}\ar@{->>}[r]& E \ar@{^{(}->}[r] & D^{\vee}}$$ where $E$ is a cocommutative $R$-bialgebra which is flat and finitely generated as an $R$-module because $D^{\vee}$ is. Consider the morphism $S_{B_f^{\vee}\ast_{A_f^{\vee}}C_f^{\vee}}:B_f^{\vee}\ast_{A_f^{\vee}}C_f^{\vee}\to B_f^{\vee}\ast_{A_f^{\vee}}C_f^{\vee}$ constructed in Remark \ref{remANTIPODO}; the commutative diagram
$$\xymatrix{B_f^{\vee}\ast_{A_f^{\vee}}C_f^{\vee}\ar[d]_{S_{B_f^{\vee}\ast_{A_f^{\vee}}C_f^{\vee}}}\ar@{->>}[r]& E \ar@{^{(}->}[r] & D^{\vee}\ar[d]^{S_{D^{\vee}}}\\ B_f^{\vee}\ast_{A_f^{\vee}}C_f^{\vee}\ar@{->>}[r]& E \ar@{^{(}->}[r] & D^{\vee} }$$
 induces a anti-morphism of $R$-algebras $$S_{E}:E\to E$$ which gives $E$ a natural structure of $R$-Hopf algebra: indeed $m_{E}\circ(id_{E} \otimes S_{E})\circ m_{E}=u_E\circ \varepsilon_E$ and $m_{E}\circ(S_{E} \otimes id_{E})\circ m_{E}=u_E\circ \varepsilon_E$ since the same equalities hold for $D^{\vee}$. It is now sufficient to take the union of $Spec(E^{\vee})$ and $N_{K}\simeq L_{K}$ in order to construct a quasi-finite and flat $R$-group scheme $P$  which is certainly a pushout in the category of quasi-finite and flat $R$-group schemes. Arguing as in the proof of Theorem \ref{teoFinito} we can deduce that $P$ is also a pushout in the category of affine $R$-group schemes. 
\endproof

\begin{remark}\label{remNonDipendeQFinito}As in Remark \ref{remNonDipendeFinito} one observes that if $U'$ is any other upper bound for $M$ and $N$ then there is a canonical isomorphism $M\sqcup_U N\simeq M\sqcup_{U'} N$: indeed $E$, as constructed in the proof, is the only quotient of $B_f^{\vee}\ast_{R}C_f^{\vee}$, $R$-flat which over $K$ gives ${B_f^{\vee}}_K\ast_{K}{C_f^{\vee}}_K\to {D^{\vee}}_K$ and this does not depend on $A_f$. The same will hold for Corollary \ref{corQFinito1} and will be used in Corollary \ref{corQFinito2}.
\end{remark}

\begin{corollary}\label{corQFinito1}
When  $N_{K}$ is \'etale then after possibly a finite extension of scalars the pushout of (\ref{eqUno}) in the category of affine $R$-group schemes exists. Again  $M\sqcup_U N$ is a lower bound for $M$ and $N$.
\end{corollary}
\proof Clear since after possibly a finite extension $K'$ of $K$ the $K$-group scheme  $N_{K}$ becomes constant then it certainly admits a finite, constant (so  flat) model over $R'$, the integral closure of $R$ in $K'$.
\endproof

Let $K'$ be a finite extension  of $K$ and $R'$ the integral closure of $R$ in $K'$ then $R'$ is a complete discrete valuation ring. Assume that $W$ is a torsion-free $R$ module of finite rank $n$ then we have the following

\begin{lemma}\label{lemTool}
If $W\otimes_RR'$ is finitely generated as an $R'$-module then  $W$ is  finitely generated as an $R$-module too (thus free).
\end{lemma}  
\proof By Theorem \ref{teoKapl} $W\simeq_{R\text{-}mod} K^{\oplus n-s}\oplus R^{\oplus s}$, where $s=dim_k(W_k)$, hence $W\otimes_R R'\simeq_{R'\text{-}mod} K'^{\oplus n-s}\oplus R'^{\oplus s}$ so if $W\otimes_RR'$ is finitely generated as an $R'$-module then $n-s=0$ and we conclude.
\endproof
This will be used in the following

\begin{corollary}\label{corQFinito2}
When $R$ is complete and $N_{K}$ is \'etale then  the pushout of (\ref{eqUno}) in the category of  affine $R$-group schemes exists. Again  $M\sqcup_U N$ is a lower bound for $M$ and $N$.
\end{corollary}
\proof Again $U_f=Spec(A_f)$, $M_f=Spec(B_f)$, $N_f=Spec(C_f)$ will denote the finite part of $U$, $M$ and $N$ respectively. Let us consider the duals $A_f^{\vee}=A^{\vee\vee}$, $B_f^{\vee}=B^{\vee\vee}$ and $C_f^{\vee}=C^{\vee\vee}$ and the commutative diagram 
\begin{equation}\label{eqUltimo}
\xymatrix{ & {A_f^{\vee}}_K\ar[dl]\ar[dr] & & & A_f^{\vee}\ar[dl]\ar[dr]\ar[lll] & \\{B_f^{\vee}}_K\ar[dr] & & {C_f^{\vee}}_K\ar[dl] & B_f^{\vee}\ar@/_/[lll] \ar[dr]& & C_f^{\vee}\ar[dl]\ar@/^/[lll]\\ & {B_f^{\vee}}_K\ast_{{A_f^{\vee}}_K}{{C_f^{\vee}}_K}\ar@{->>}[d] & & & B_f^{\vee}\ast_{A_f^{\vee}}C_f^{\vee} \ar[lll]& \\  & E\ar@{^{(}->}[d] & & & & \\ & C^{\vee} & & & & }
\end{equation}
where $E$ comes from the factorisation of the universal morphism ${B_f^{\vee}}_K\ast_{{A_f^{\vee}}_K}{{C_f^{\vee}}_K}\to C^{\vee}$. Arguing as in Theorem \ref{teoQFinito} we provide it with a natural structure of $K$-Hopf algebra. Using again  \cite{EGAIV2}, Lemme (2.8.1.1) we construct the unique quotient $${B_f^{\vee}}\ast_{{A_f^{\vee}}}{{C_f^{\vee}}}\twoheadrightarrow E'$$ which is $R$-flat and which generically  gives $${B_f^{\vee}}_K\ast_{{A_f^{\vee}}_K}{{C_f^{\vee}}_K}\twoheadrightarrow E.$$
Thus $E'$ has naturally a structure of a cocommutative $R$-Hopf algebra: indeed it inherits from ${B_f^{\vee}}\ast_{{A_f^{\vee}}}{{C_f^{\vee}}}$ a cocommutative $R$-coalgebra structure and by means of \cite{EGAIV2}, (2.8.3) an anti-morphism of $R$-agebras $S_E':E'\to E'$ which is a coinverse since tensoring it over $K$ we obtain $S_E:E\to E$ which is a coinverse for $E$. If we prove that $E'$ is finitely generated as a $R$-module then $Spec(E'^{\vee})$ glued to $N$ is the desired pushout. So now it remains to prove that $E'$ is finitely generated as a $R$-module: let $K\to K'$ be a finite field extension such that $N_{K'}$ admits a finite and flat model over  $R'$, the integral closure of $R$ in $K'$. Then  by Corollary \ref{corQFinito1}   and Remark \ref{remNonDipendeQFinito} $E'\otimes_RR'$ is $R$-finite and flat. Lemma \ref{lemTool} implies that $E'$ is $R$-finite and flat too.
\endproof

\begin{remark}\label{remLowBoundFinito}
It is less elegant but still true that Corollary \ref{corQFinito2} holds for all those $N_K$ that admits a finite and flat $R$-model after possibly a finite extension of scalrs and \'etale ones are just a particular case. Observe furthermore that, following the proof, in the situation of both Theorem \ref{teoQFinito} and Corollary \ref{corQFinito1} one can find a finite and flat lower bound for $M$ and $N$. This can be false in the situation of Corollary \ref{corQFinito2}. 
\end{remark}

\subsection{Cokernels and quotients}
\label{sez:cok}
In a category $\mathcal{C}$ with zero object  $0_{\mathcal{C}}$ (that is an object wich is both initial and final), we can define  the cokernel of a morphism $f:A\to B$ (see for instance \cite{ML}, III, \S 3) which turns out to be the pushout $0_{\mathcal{C}}\sqcup_A B$ of the obvious diagram. As explained in the introduction in this section we are going to describe, in Proposition \ref{propEsisteCok}, a new and easy proof for a well known result. First we need a lemma:

\begin{lemma}\label{lemStarProd}
Let $R$ be a Dedekind ring or a field, $A$, $B$ and $C$ $R$-Hopf algebras provided with $R$-Hopf algebra morphisms $f:A\to B$ and $g:A\to C$. Then the star product $B\star_AC$ defined in Definition \ref{defQuozTensor} has a natural structure of $R$-Hopf algebra.
\end{lemma}
\proof First we prove the existence of the comultiplication $\Delta_{B\star_AC}$: it is sufficient to consider the following diagram
$$\xymatrix{& B\ar[rd]^u\ar[rrr]^{\Delta_B} & & & B\otimes_R B \ar[rd]^{u\otimes u} & \\A\ar[ru]^f\ar[rd]_g\ar@/^/[rrr]^(.3){\Delta_A} & & B\star_AC\ar@{-->}@/_/[rrr]_(.6){\Delta_{B\star_AC}} & A\otimes_R A \ar[ru]^{f\otimes f}\ar[rd]_{g\otimes g}& & (B\star_AC)\otimes_R (B\star_AC) \\ & C\ar[ru]_v\ar[rrr]_{\Delta_C} & & & C\otimes_R C \ar[ru]_{v\otimes v} & } $$
where the existence of $\Delta_{B\star_AC}$ is ensured by Proposition \ref{propUnivStar}. The existence of $\varepsilon_{B\star_AC}$ is easier and an argument similar to the one used in Remark \ref{remANTIPODO} ensures the existence of an anti-morphism $S_{B\star_AC}:{B\star_AC}\to {B\star_AC}$ which is compatible with $S_{B\otimes_RC}$, i.e. if $\lambda:{B\otimes_RC}\twoheadrightarrow {B\star_AC}$ denotes the canonical projection then $\lambda \circ S_{B\otimes_RC}=S_{B\star_AC}\circ \lambda$. From this we deduce that $S_{B\star_AC}$ is the desired coinverse for ${B\star_AC}$.
%This implies that $$m_{B\star_AC}\circ (S_{B\star_AC}\otimes id_{B\star_AC})\circ \Delta_{B\star_AC}\circ \lambda =\lambda \circ m_{B\otimes_RC}\circ (S_{B\otimes_RC}\otimes id_{B\otimes_RC})\circ \Delta_{B\otimes_RC}$$ which, in turn, implies that $m_{B\star_AC}\circ (S_{B\star_AC}\otimes id_{B\star_AC})\circ \Delta_{B\star_AC}$ is a $R$-algebra morphism since $m_{B\otimes_RC}\circ (S_{B\otimes_RC}\otimes id_{B\otimes_RC})\circ \Delta_{B\otimes_RC}$ is. Then it is the unique morphism making the following diagram
\endproof

\begin{proposition}\label{propEsisteCok}Let $R$ be a Dedekind ring or a field, $G$ and $H$ two finite and flat $R$-group schemes and $f:H\to G$ a morphism of $R$-group schemes. Then the cokernel of $f$ exists in the category of $R$-affine group schemes.
\end{proposition}
\proof The zero object in the category of $R$-affine group schemes is $Spec(R)$. Let us set $H=Spec(A)$ and $G=Spec(B)$. Then we first compute the pushout in the category of $R$-Hopf algebras of the diagram 
\begin{equation}\label{eqPusherHopf}
\xymatrix{ & A^{\vee}\ar[dl]\ar[dr] & \\ B^{\vee} & & R .}
\end{equation}
In Example \ref{exaCopro2}  we have observed that $W:=B^{\vee}\ast_{A^{\vee}}R\simeq B^{\vee}\star_{A^{\vee}}R$ canonically. That $W$ has a natural structure of $R$-Hopf algebra follows from Lemma \ref{lemStarProd}. If $R$ is a Dedekind ring and  $W$ is not flat then we consider $F(W)$ (cf. Notation \ref{notFM}) which is flat and finitely generated and inherits the $R$-Hopf algebra structure. Sice the case of a field is similar let us just consider the case of a Dedekind ring $R$: we are now going to prove that $Spec (F(W)^{\vee})$ is the desired pushout. So let $M$ be any affine $R$-group scheme, $v:G\to M$ an $R$-group scheme morphism and $u:Spec(R)\to M$ the natural inclusion (the unity map). Let us assume we have a commutative diagram
\begin{equation}\label{eqPusher}
\xymatrix{ & Spec(R) \ar[dr]^u & \\H\ar[ur]\ar[dr]_f & & M \\ & G\ar[ur]_v & }
\end{equation}
Observe that we can assume $M$ to be finite and flat, for if it is not we can factor $v$ through a  finite and flat (since $G$ is)  $R$-group scheme that makes a diagram similar to (\ref{eqPusher}) commmute. When $M$ is finite and flat it is easy to construct a universal morphism $Spec(F(W)^{\vee})\to M$ since $F(W)$ is easily seen to be the pushout  of diagram (\ref{eqPusherHopf}) in $R$-$\mathcal{H}opf_{ff}$.
\endproof

\begin{corollary}\label{corEsisteCok}
Let $R$ be a Dedekind ring, $G$ and $H$ two finite and flat $R$-group schemes with $H$ a closed and normal $R$-subgroup scheme of $G$. Then the quotient $G/H$ exists in the category of $R$-affine group schemes.
\end{corollary}
\proof This follows directly from Proposition \ref{propEsisteCok} where we take for $f:H\to G$ the given closed immersion.
\endproof
%{\color{red} A lower bound } for models of finite group schemes...che poi alla fine e` quello che mi interessa... .{\color{red} cokernel} e poi ancora {\color{red} quoziente} 

\begin{flushright}Marco Antei\\ 

Department of Mathematics \\
Ben Gurion University of the Negev\\
Be'er Sheva 84105, Israel\\
E-mail:  \texttt{anteim@math.bgu.ac.il}\\
\texttt{marco.antei@gmail.com}
\end{flushright}
\end{document}